\documentclass[a4paper, 11pt, ]{amsart}

\newcommand{\papertitle}{An Induction Property for Prime Counting Functions}

\usepackage{times,amsmath,amsthm,amssymb,color,enumerate,accents,csquotes,xypic,placeins,mathtools}
\usepackage[colorlinks=true,hyperindex,pagebackref=false, citecolor=cyan]{hyperref}
\usepackage[T1]{fontenc}
\usepackage[utf8]{inputenc}
\usepackage[russian,english]{babel}

\theoremstyle{remark}
\newtheorem*{remark*}{Remark}

\newcommand{\QQ}{\mathbb Q}

\renewcommand{\c}{\mathfrak c}

\newcommand{\Gal}{\operatorname{Gal}}

\newcommand{\nf}{K}
\newcommand{\enf}{L}
\newcommand{\mnf}{E}

\newcommand{\mfp}{\mathfrak p}

\newcommand{\mfP}{\mathfrak P}

\newcommand{\N}{\mathcal N}

\newcommand{\nfps}[1]{\Sigma_{#1}} 
\newcommand{\s}{S} 

\newcommand{\ind}{\operatorname{Ind}} 

\newcommand{\tik}{\[\begin{tikzcd}}
\newcommand{\zcd}{\end{tikzcd}\]}

\newcommand{\chebotarev}{Chebotar\"ev}

\newcommand{\Li}{\operatorname{Li}}
\newcommand{\fbnf}{\psi} 
\newcommand{\clfn}{\varphi} 
\newcommand{\E}[2]{{\mathbb E}_{#1}[{#2}]} 

\newcommand{\indprop}{induction property } 
\newcommand{\Ch}{\foreignlanguage{russian}{Ч}}

\makeatletter
\newcommand{\oset}[3][0ex]{%
  \mathrel{\mathop{#3}\limits^{
    \vbox to#1{\kern-2\ex@
    \hbox{$\scriptstyle#2$}\vss}} } }
\makeatother

\newcommand{\grant}{The contents of this paper were worked out while the author was partially supported by the National Science Foundation via grant DMS-1601844.}

\title{\papertitle}
\author{Andrew O'Desky}
\date{\today}

\begin{document}

\begin{abstract}
    We provide an elementary proof of an asymptotic formula for prime counting functions. As a minor application we give a new reduction of the proof of {\chebotarev}'s density theorem to the cyclic case. 
\end{abstract}

\maketitle


Let $\nf$ be a number field, $\nfps{\nf}$ the set of its finite primes, and $\Omega$ a set. 
A function $\fbnf\colon \nfps{\nf} \backslash \s \to \Omega$ defined away from a finite set $\s$ is \textit{frobenian} if there exist a finite Galois extension $\enf/\nf$ unramified outside $\s$ and an $\Omega$-valued class function $\clfn$ of $G:=\Gal(\enf/\nf)$ such that $\fbnf(\mfp) = \clfn(\sigma_\mfp)$ for all unramified $\mfp$ where $\sigma_\mfp$ is the conjugacy class of Frobenius elements of primes of $\enf$ dividing $\mfp$. 
This notion is due to Serre \cite[\S 3.3]{serre}.
Examples of frobenian functions (for varying $\Omega$) include $p \mod N$ for a positive integer $N$, the number of roots mod $p$ of a non-zero integral polynomial in one variable, and the $p$-th Fourier coefficient of a modular form. 
Frobenian functions are closely related to the notions of Frobenius sets and {\chebotarev} sets for which we refer the reader to \cite{lagarias}, \cite{cheb}, \cite{rosen}. 

The aim of this note is to demonstrate an elementary proof of the following formula for real-valued frobenian functions: 
\begin{equation} \label{eqn:1} \tag{$\star$}
    \E{\mnf}{\fbnf}
=
\E{\nf}{\fbnf'}
\end{equation}
where $\mnf$ is an intermediate extension of $\enf/\nf$, $\fbnf$ is a frobenian function for $\mnf$, $\fbnf'$ is the induction\footnote{the frobenian function for $\nf$ whose associated class function is $\ind^G_H(\clfn)$, $H=\Gal(\enf/\mnf)$.} of $\fbnf$ to $K$, 
$$
\E{\mnf}{\fbnf}
:=
\lim_{x \to \infty}
\frac
    {\sum_{\N\mfp \leq x} \fbnf(\mfp) }
    {\int_2^x (\log t)^{-1}\, dt}, 
$$ 
and the sum is over the primes $\mfp \in \nfps{\mnf}\backslash \s$ such that $\N \mfp \leq x$. 
Formulas for $\E{\mnf}{\fbnf}$ imply various classical density theorems as special cases, such as Dirichlet's theorem on primes in arithmetic progression, Landau's prime ideal theorem, and {\chebotarev}'s density theorem. 
We will refer to \eqref{eqn:1} as the induction property for prime counting functions. 
A quick proof of \eqref{eqn:1} proceeds as follows: 
rearrange the sum over primes according to their conjugacy class of Frobenius elements and use {\chebotarev}'s density theorem to see that $\E{\mnf}{\fbnf}$ and $\E{\nf}{\fbnf'}$ have the limiting values 
$\sum_{C \in \c(H)} 
    \clfn(C) 
        {|C|}
        {|H|} ^{-1}
$ and
$
\sum_{C \in \c(G)} 
    \clfn'(C) 
        {|C|}
        {|G|}^{-1} 
        $ where $H:= \Gal(\enf/\mnf)$ and $\clfn' := \ind^G_H(\clfn)$.
These quantities are equal by Frobenius reciprocity.


This note demonstrates that \eqref{eqn:1} admits an elementary proof which does not use {\chebotarev}'s density theorem or Frobenius reciprocity, which is somewhat surprising in consideration of the proof just given.
Bypassing {\chebotarev}'s theorem will require a more intricate argument using other tools from algebraic number theory, namely Chebyshev's bounds on $\pi(x)$ and the double-coset description of prime factorization in an intermediate field extension. 
The key idea in the proof is relating the Frobenius elements in the support of an induced class function to primes with inertia degree one (see \eqref{eqn:ind}). 
Afterward, we turn around in \S\ref{sec:reduction} and use our proof of \eqref{eqn:1} to give a minor application by reducing the proof of {\chebotarev}'s theorem to the cyclic case. 
Deuring also gave a cyclic reduction of {\chebotarev}'s theorem \cite{deuring} (and see \cite{maccluer} for a similar proof to Deuring's). 
The \indprop \eqref{eqn:1} may be interpreted as the ``underlying mechanism'' which makes the reduction possible.

\section{Elementary Proof of \eqref{eqn:1}} 

As the portion of the sum $\sum_{\N\mfp \leq x} \fbnf(x)$ corresponding to primes of $\mnf$ with inertia degree $>1$ over $\QQ$ is $O(\pi(x^{1/2}))$, any modestly effective asymptotic estimate for $\pi(x)$ shows that such primes will not contribute in the limit.\footnote{For instance, Chebyshev's bound that $\frac{Ax}{\log x} \leq \pi(x) \leq \frac{Bx}{\log x}$ for some $0 < A < B$, $x \gg 0$. Note that $\Li(x)x^{-1}\log x \to 1$ as $x \to \infty$.}  
 Let $\nfps{\mnf}^1 \subset \nfps{\mnf}$ be the subset of primes of $\mnf$ with inertia degree one over $\QQ$. 
 Then we have equality of $\E{\mnf}{\fbnf}$ with  
 \begin{equation*}
\lim_{x \to \infty} 
    \sum_{p \leq x} 
    \sum_{\substack{\mfP \in \nfps{\mnf}^1,\\ \mfP / p}} 
        \fbnf(\mfP) 
        \cdot 
        \Li(x)^{-1}, 
 \end{equation*}
 where we omit the finitely many primes for which $\fbnf$ is undefined for the sake of simplicity. 
 Grouping primes $\mfP$ lying over the same prime $\mfp$ of $K$ yields the equality 
 \begin{equation}\label{eqn:sum}
     \E{\mnf}{\fbnf}
=
\lim_{x \to \infty} 
\sum_{p \leq x} 
\sum_{\substack{\mfp \in \nfps{\nf}^1,\\ \mfp / p}} 
\sum_{\substack{\mfP \in \nfps{\mnf}^1,\\ \mfP / \mfp}} 
    \fbnf(\mfP) 
    \cdot 
    \Li(x)^{-1}. 
 \end{equation}
 Now we use the double-coset description of primes in an intermediate extension to rewrite \eqref{eqn:sum} in terms of the class function $\clfn$ of $\fbnf$; see e.g. \cite[\S 1.9]{neukirch}.
 Let $\mfp$ be a prime of $\nf$ unramified in $\enf$, $\mfP$ any prime of $\enf$ dividing $\mfp$, and $s_\mfp$ the Frobenius element of $\mfP$ in $\enf/\nf$, i.e. the canonical generator of the decomposition subgroup of $\mfP$ in $G$. 
Recall that primes of $\mnf$ dividing $\mfp$ are in bijection with the double cosets in $ \langle s_\mfp \rangle \backslash G/H$ under 
$\mfP^g \cap \mnf 
\leftrightarrow 
\langle s_\mfp \rangle g H $ 
and the inertia degree of $\mfP^g \cap \mnf$ over $\mfp$ is $| \langle s_\mfp \rangle g H ||H|^{-1}$. 
Let $g_1,\ldots,g_r$ be a set of representatives for $G/H$ and set $\mfP_i = \mfP^{g_i}\cap \mnf$. 
Then 
$$
    g_i^{-1}s_\mfp g_i \in H 
\iff 
    g_i H 
    \text{ fixed by } 
    s_\mfp 
\iff 
    |\langle s_\mfp \rangle g_i H| 
    = 
    |H|
\iff
    f(\mfP_i/\mfp) = 1,
$$ 
and for such $\mfP_i$ one has that $g_i^{-1}s_\mfp g_i \in \sigma_{\mfP_i} \subset H$ (e.g. \cite[\S 1.8, Proposition 23(a)]{serre-lf}). 
Let $\clfn$ be the class function of $H$ corresponding to $\fbnf$.
We extend $\clfn$ by zero to a function on $G$ that is not necessarily constant on conjugacy classes of $G$. 
The induction of $\clfn$ is the class function $\clfn'$ of $G$ given by 
$
\clfn'(g) 
= 
\sum_{i = 1}^r
\clfn(g_i^{-1} g g_i).
$
If $f(\mfP_i/\mfp) = 1$ then $\psi(\mfP_i) = \clfn(g_i^{-1}s_\mfp g_i)$ since $g_i^{-1} s_\mfp g_i \in \sigma_{\mfP_i}$, so we obtain a formula for $\fbnf'$
\begin{equation}\label{eqn:ind}
    \fbnf'(\mfp)
    =
    \clfn'(s_\mfp)
    =
    \sum_{i = 1}^r
        \clfn(g_i^{-1} s_\mfp g_i)
    =
    \sum_{f(\mfP/\mfp) = 1}
    \fbnf(\mfP), \quad (\mfP / \mfp,\, \mfP \in \nfps{\mnf}).
\end{equation}
Note that the $\mfP_i$ are not necessarily distinct but the $\mfP_i$ with inertia degree one over $\mfp$ are distinct. 
Combining \eqref{eqn:sum} and \eqref{eqn:ind},
\begin{align*}
\sum_{p \leq x} 
\sum_{\substack{\mfp \in \nfps{\nf}^1,\\ \mfp / p}} 
&\sum_{\substack{\mfP \in \nfps{\mnf}^1,\\ \mfP / \mfp}} 
    \fbnf(\mfP) 
    \cdot 
    \Li(x)^{-1} 
=  
\sum_{p \leq x} 
\sum_{\substack{\mfp \in \nfps{\nf}^1,\\ \mfp / p}} 
\sum_{i = 1}^r 
\clfn(g_i^{-1}s_\mfp{g_i}) 
    \cdot 
    \Li(x)^{-1} \\
&= \sum_{p \leq x} 
\sum_{\substack{\mfp \in \nfps{\nf}^1,\\ \mfp / p}} 
   \clfn'(s_\mfp)
   \cdot 
   \Li(x)^{-1} \\
&= \sum_{p \leq x} 
\sum_{\substack{\mfp \in \nfps{\nf}^1,\\ \mfp / p}} 
   \fbnf'(\mfp) 
   \cdot 
   \Li(x)^{-1} \qquad \text{(definition of $\fbnf'$)}. 
\end{align*}
By the same technique as in the beginning of the proof we may add in the primes of $K$ with inertia degree $>1$ without changing the limiting value of the above sums and so it converges to $\E{\nf}{\fbnf'}$ as $x \to \infty$.
Combining this with \eqref{eqn:sum} finishes the proof of \eqref{eqn:1}.

\begin{remark*}
    The key point of the proof is \eqref{eqn:ind} which relates the support of the induced class function $\clfn'$ with the contribution to $\E{\nf}{\fbnf'}$ from primes with inertia degree one over $\QQ$. 
    In some sense this shows that induction of class functions is ``built into'' prime factorization in number fields.
\end{remark*}

\begin{remark*}
There is also an analytic version -- in the sense of analytic density, 
\begin{equation} \label{eqn:2} \tag{$\star\star$}
\lim_{s \to 1^+} 
\sum_{\mfP \in \nfps{\mnf}\backslash \s'} 
        \frac
            {\fbnf(\mfP)}
            {\N(\mfP)^{s}}
        \cdot z(s)^{-1} 
= 
\lim_{s \to 1^+} 
\sum_{\mfp \in \nfps{\nf}\backslash \s} 
    \frac
        {\fbnf'(\mfp)}
        {\N (\mfp)^{s}} 
    \cdot z(s)^{-1}  
\end{equation}
where $z(s) = -\log(s-1)$.
As one might expect the proof of \eqref{eqn:2} is analogous to the proof above so we omit it. The key difference is that in place of Chebyshev's bound one uses the fact that $\sum p^{-f}$ only diverges when $f = 1$ and converges absolutely otherwise.
\end{remark*}


\section{Proof of Reduction for {\chebotarev}'s density theorem}\label{sec:reduction} 

\emph{Notation.} 
As above let $\enf/\nf$ be a finite Galois extension with Galois group $G$ unramified outside a finite set $\s \subset \nfps{\nf}$, and let $\mnf$ be an intermediate extension with $H:=\Gal(L/E)$. 
For any conjugacy class $C \subset G$ let $\Gamma^\nf_C(x) = |\{ \mfp \in \nfps{\nf} : \sigma_\mfp = C, \N(\mfp) \leq x \}|$ and let $d(\Gamma^\nf_C) := \lim_{x \to \infty} \Gamma_C^K(x) \Li(x)^{-1}$.\\[-3mm]

Assume that {\chebotarev}'s density theorem is true for all cyclic extensions for varying $K$. 
We will show that $d(\Gamma^\nf_C) = |C||G|^{-1}$.
We apply {\chebotarev}'s theorem in the form 
\begin{equation}\label{eqn:CD}\tag{\Ch}
    \E{\mnf}{\fbnf}
= 
|H|^{-1} 
\sum_{h \in H} 
    \clfn(h)
\end{equation}
for any frobenian function $\fbnf$ for $\mnf$ whenever $H$ is cyclic, where $\clfn$ is the class function associated with $\fbnf$. 

The proof is by induction on the order $n$ of the elements in $C$. 
Let $\sigma\in C$, $E$ the fixed field of $\sigma$, and $H = \langle \sigma \rangle$. 
For an integer $e$ let $c_e$ be the order of the centralizer of $\sigma^e$ in $G$. 
Consider the class function of $H$ given by 
$$
\clfn 
= 
\sum_{e|n} 
\delta_{\sigma^e}
nc_e^{-1}
$$ 
where $\delta_{\sigma^e}$ is a delta function supported at $\sigma^e$; 
then the induction of $\clfn$ to $G$ equals $\sum_{e | n} {1}_{C^e}$ where $C^e$ is the conjugacy class of $\sigma^e$. 
When $\sigma = 1$ applying the \indprop \eqref{eqn:1} and \eqref{eqn:CD} for $E = L$ yields 
$\E{\enf}{\fbnf} = |G|^{-1} = \E{\nf}{\fbnf'} = d(\Gamma_{\{1\}}^K),$ 
thereby proving {\chebotarev}'s density theorem in this case. 
Now assume $\sigma \neq 1$. 
We have that 
\begin{align*}  
    \E{\mnf}{\fbnf} &=|H|^{-1} \sum_{h \in H} \clfn(h) && \eqref{eqn:CD} \\ 
&= |G|^{-1}\sum_{g \in G} \clfn'(g) && \text{(Frobenius reciprocity)} \\
&= \sum_{e|n} \frac{|C^e|}{|G|} &&\text{(above formula for $\clfn'$).} \\
\end{align*} 
By the inductive hypothesis $\E{\nf}{\fbnf'} = d(\Gamma_C^K) + \sum_{e | n, e >1}{|C^e|}{|G|}^{-1}$. 
By the inductive property \eqref{eqn:1} this equals $\E{\mnf}{\fbnf} = \sum_{e|n} {|C^e|}{|G|}^{-1}$. 
Taking their difference shows that $d(\Gamma^K_C) = |C| |G|^{-1}$.

\begin{remark*}
A reduction to the cyclic case of the weaker analytic version of {\chebotarev}'s theorem may also be deduced along the same lines with \eqref{eqn:2} used in place of \eqref{eqn:1}.
\end{remark*}


   \section{Acknowledgments} 
   The author thanks Julian Rosen and Michael Zieve for pointing out the sufficiency of a weaker asymptotic for $\pi(x)$ in the proof of \eqref{eqn:1}. 
{\grant}


\medskip

\bibliographystyle{plain}
\bibliography{induction}

\end{document}